\documentclass[twoside,11pt]{article}
\usepackage{graphicx,amsmath,amsthm,amssymb,latexsym,amsfonts,color,cite}
\usepackage[bookmarksnumbered, colorlinks, plainpages]{hyperref}
\usepackage{cases}

\newtheorem{thm}{\bf Theorem}

\newtheorem{example}{\bf Example}

\newtheorem{remark}{\bf Remark}

\begin{document}
\title{\bf A tighter $S$-type singular value inclusion set for rectangular tensors}
\author{{Caili Sang\footnote{Corresponding author, E-mail: sangcl@126.com; sangcaili@gzmu.edu.cn (Caili Sang)}}\\ [2mm]
{\small \textit{College of Data Science and Information Engineering, Guizhou Minzu University,}}\\
{\small \textit{Guiyang 550025, P.R.China}}\\
}

\date{}
\maketitle
\noindent {\bf Abstract.}
A new $S$-type singular value inclusion set for rectangular tensors is given and proved to be tighter than that in [Sang C.L., An $S$-type singular value inclusion set for rectangular tensors, J. Inequal. Appl. 2017: 141, 2017]. Based on this set, new bounds for the largest singular value of nonnegative rectangular tensors are obtained and proved to be better than some existing results.
Compared with the results in the paper mentioned above, the advantage of the new results is that, under the same computations, we can obtain a tighter singular value inclusion set for rectangular tensors and sharper bounds for the largest singular value of nonnegative rectangular tensors. Finally,
a numerical example is given to verify the theoretical results.  \\[-3mm]

\noindent{\it\bf  Keywords}: {rectangular tensor; nonnegative tensors; singular value; inclusion set}\\
\noindent{\it\bf  AMS Subject Classification}: 15A18; 15A42; 15A69. \\
\pagestyle{myheadings}\markboth{Caili Sang}{A tighter $S$-type Singular value inclusion sets for rectangular tensors}
\thispagestyle{empty}

\section{Introduction}
Let $\mathbb{R} (\mathbb{C})$ be the real (complex) field, $p,q,m,n$ be positive integers, $l=p+q$, $m,n\geq 2$ and $N=\{1,2,\cdots,n\}$.
$\mathcal{A}=(a_{i_1\cdots i_pj_1\cdots j_q})$ is called a real $(p,q)$th order $m\times n$ dimensional rectangular tensor, or simply a real rectangular tensor, denoted by $\mathcal{A}\in\mathbb{R}^{[p,q;m,n]}$, if
\[a_{i_1\cdots i_pj_1\cdots j_q}\in \mathbb{R},~1\leq i_1\cdots i_p\leq m,~1\leq j_1\cdots j_q\leq n. \]
When $p=q=1, \mathcal{A}$ is simply a real $m\times n$ rectangular matrix.
$\mathcal{A}$ is called a nonnegative rectangular tensor, denoted by $\mathcal{A}\in\mathbb{R}^{[p,q;m,n]}_{+}$, if each of its entries
$a_{i_1\cdots i_pj_1\cdots j_q}\geq 0$.

If there are
a number $\lambda\in\mathbb{C},$ vectors $x=(x_1,x_2,\cdots, x_m)^\textmd{T}\in\mathbb{C}^m\backslash\{0\}$ and $y=(y_1,y_2,\cdots, y_n)^\textmd{T}\in \mathbb{C}^n\backslash\{0\}$ such that
\begin{equation*}
\left\{
\begin{aligned}
\mathcal{A}x^{p-1}y^q=\lambda x^{[l-1]},\\
\mathcal{A}x^py^{q-1}=\lambda y^{[l-1]},
\end{aligned}
\right.
\end{equation*}
then $\lambda$ is called a singular value of $\mathcal{A}$, and $(x,y)$ is the left and right eigenvectors
pair of $\mathcal{A}$, associated with $\lambda$, respectively, where
$\mathcal{A}x^{p-1}y^q$ and $x^{[l-1]}$ are vectors in $\mathbb{R}^m$, whose $i$th component is
$$(\mathcal{A}x^{p-1}y^q)_i=\sum\limits_{i_2,\cdots,i_p=1}^m\sum\limits_{j_1,\cdots,j_q=1}^na_{ii_2\cdots i_pj_1\cdots j_q}x_{i_2}\cdots x_{i_p}y_{j_1}\cdots y_{j_q}$$
and
$$(x^{[l-1]})_i=x_i^{l-1};$$
$\mathcal{A}x^py^{q-1}$ and $x^{[l-1]}$ are vectors in $\mathbb{R}^n$ whose $j$th component is
$$(\mathcal{A}x^py^{q-1})_j=\sum\limits_{i_1,\cdots,i_p=1}^m\sum\limits_{j_2,\cdots,j_q=1}^na_{i_1\cdots i_pjj_2\cdots j_q}x_{i_1}\cdots x_{i_p}y_{j_2}\cdots y_{j_q}$$
and
$$(y^{[l-1]})_j=y_j^{l-1}.$$
Furthermore, if $\lambda\in \mathbb{R}, x\in \mathbb{R}^m,$ and
$y\in\mathbb{R}^n$, then we say that $\lambda$ is an H-singular value of $\mathcal{A}$, and $(x,y)$ is the left
and right H-eigenvectors pair associated with $\lambda$, respectively. If a singular value
is not an H-singular value, we call it an N-singular value of $\mathcal{A}$.
We call
$$\lambda_0=\max\{|\lambda|:\lambda~is~ a~ singular~ value~ of~ \mathcal{A}\}$$
is the largest singular value; see \cite{changkc,yangqz,lim} for details.

When $p=q,m=n$, such real rectangular tensors have a wide range of practical applications in
the strong ellipticity condition problem in solid mechanics \cite{jk,yw} and the entanglement problem in quantum physics
\cite{dd,einstein}. For example, the elasticity tensor is a tensor with $p=q=2$ and $m=n=2$ or 3; for details, see \cite{changkc}.


Because it is not easy to compute all singular values (eigenvalues) of tensors when the order
and dimension are large, one always tries to give a set including all singular values (eigenvalues) in the complex plane \cite{lcq-jaq,lcq-lyt,lcq,scljia,zjxjia}, or give upper and lower bounds for the largest singular value of nonnegative rectangular tensors\cite{hj,zjxjias,zjx}.
Very recently, Sang \cite{scljia} proposed the following $S$-type singular value inclusion set for rectangular tensors by breaking $N=\{1,2,\cdots,n\}$ into disjoint subsets $S$ and
its complement $\bar{S}$.
\begin{thm}\emph{\cite[Theorem 1]{scljia}}\label{scl-th1}
Let $\mathcal{A}\in \mathbb{R}^{[p,q;n,n]}$, $S$ be a nonempty proper subset of $N$, $\bar{S}$ be the complement of $S$ in N. Then
$$\sigma(\mathcal{A})\subseteq \Upsilon^S(\mathcal{A})=
\left(\bigcup\limits_{i\in S,j\in \bar{S}}\left(\hat{\Upsilon}_{i,j}(\mathcal{A})\bigcup \tilde{\Upsilon}_{i,j}(\mathcal{A})\right)\right)
\bigcup
\left(\bigcup\limits_{i\in \bar{S},j\in S}\left(\hat{\Upsilon}_{i,j}(\mathcal{A})\bigcup \tilde{\Upsilon}_{i,j}(\mathcal{A})\right)\right),
$$
where
\begin{eqnarray*}
\hat{\Upsilon}_{i,j}(\mathcal{A})=\left\{z\in \mathbb{C}: (|z|-r_i^j(\mathcal{A}))|z|\leq |a_{ij\cdots jj\cdots j}|\max\{R_j(\mathcal{A}),C_j(\mathcal{A})\}\right\},\\
\tilde{\Upsilon}_{i,j}(\mathcal{A})=\left\{z\in \mathbb{C}: (|z|-c_i^j(\mathcal{A}))|z|\leq |a_{j\cdots jij\cdots j}|\max\{R_j(\mathcal{A}),C_j(\mathcal{A})\}\right\},
\end{eqnarray*}
and
\begin{eqnarray*}
&&R_i(\mathcal{A})=\sum\limits_{i_2,\cdots,i_p,j_1,\cdots,j_q\in N}|a_{ii_2\cdots i_pj_1\cdots j_q}|,\\
&&r_i^j(\mathcal{A})=\sum\limits_{\delta_{ji_2\cdots i_pj_1\cdots j_q}=0}|a_{ii_2\cdots i_pj_1\cdots j_q}|=R_i(\mathcal{A})-|a_{ij\cdots jj\cdots j}|,\\
&&C_j(\mathcal{A})=\sum\limits_{i_1,\cdots,i_p,j_2,\cdots,j_q\in N}|a_{i_1\cdots i_pjj_2\cdots j_q}|,\\
&&c_j^i(\mathcal{A})=\sum\limits_{\delta_{i_1\cdots i_pij_2\cdots j_q}=0}|a_{i_1\cdots i_pjj_2\cdots j_q}|=C_j(\mathcal{A})-|a_{i\cdots iji\cdots i}|,\\
&&\delta_{i_1\cdots i_pj_1\cdots j_q}=\left\{\begin{array}{rl}1&if~i_1=\cdots=i_p=j_1=\cdots=j_q,\\0&otherwise.\end{array}\right.
\end{eqnarray*}
\end{thm}

Based on the set in Theorem \ref{scl-th1}, Sang in \cite{scljia} obtained the following upper and lower bounds for the largest singular value of nonnegative rectangular tensors.
\begin{thm}\emph{\cite[Theorem 2]{scljia}}\label{scl-th2}
Let $\mathcal{A}=(a_{i_{1}\cdots i_{m}})\in \mathbb{R}^{[p,q;n,n]}_{+}$, $S$ be a nonempty proper subset of $N$, $\bar{S}$ be the complement of $S$ in $N$. Then
\begin{eqnarray*}
L^S(\mathcal{A})\leq\lambda_0\leq U^S(\mathcal{A}),
\end{eqnarray*}
where
\begin{eqnarray*}
L^S(\mathcal{A})=\min\{\hat{L}^S(\mathcal{A}),\hat{L}^{\bar{S}}(\mathcal{A}),\tilde{L}^S(\mathcal{A}),\tilde{L}^{\bar{S}}(\mathcal{A})\}\\ U^S(\mathcal{A})=\max\{\hat{U}^S(\mathcal{A}),\hat{U}^{\bar{S}}(\mathcal{A}),\tilde{U}^S(\mathcal{A}),\tilde{U}^{\bar{S}}(\mathcal{A})\},
\end{eqnarray*}
and
\begin{eqnarray*}
&&\hat{L}^S(\mathcal{A})=\min\limits_{i\in S,j\in \bar{S}}\frac{1}{2}\left\{r_i^j(\mathcal{A})
+[(r_i^j(\mathcal{A}))^2+4a_{ij\cdots jj\cdots j}\min\{R_j(\mathcal{A}),C_j(\mathcal{A})\}]^{\frac{1}{2}}\right\},\\
&&\tilde{L}^S(\mathcal{A})=\min\limits_{i\in S,j\in \bar{S}}\frac{1}{2}\left\{c_i^j(\mathcal{A})
+[(c_i^j(\mathcal{A}))^2+4a_{j\cdots jij\cdots j}\min\{R_j(\mathcal{A}),C_j(\mathcal{A})\}]^{\frac{1}{2}}\right\},\\
&&\hat{U}^S(\mathcal{A})=\max\limits_{i\in S,j\in \bar{S}}\frac{1}{2}\left\{r_i^j(\mathcal{A})
+[(r_i^j(\mathcal{A}))^2+4a_{ij\cdots jj\cdots j}\max\{R_j(\mathcal{A}),C_j(\mathcal{A})\}]^{\frac{1}{2}}\right\},\\
&&\tilde{U}^S(\mathcal{A})=\max\limits_{i\in S,j\in \bar{S}}\frac{1}{2}\left\{c_i^j(\mathcal{A})
+[(c_i^j(\mathcal{A}))^2+4a_{j\cdots jij\cdots j}\max\{R_j(\mathcal{A}),C_j(\mathcal{A})\}]^{\frac{1}{2}}\right\}.
\end{eqnarray*}
\end{thm}

In this paper, by the technique in \cite{lcq-jaq}, we give a new $S$-type singular value inclusion set for a real rectangular tensor $\mathcal{A}$ and prove that the new set is tighter than $\Upsilon^S(\mathcal{A})$. As an application, we obtain new upper and lower bounds for the largest singular value of nonnegative rectangular tensors and prove that the new  bounds are better than those in Theorem \ref{scl-th2} and Theorem 4 of \cite{yangqz}.
\section{Main results}\label{sec2}

We begin with some notation. Given a nonempty proper subset $S$ of $N$, we denote
\begin{eqnarray*}
&&\Delta^N:=\{(i_2,\cdots i_p,j_1,\cdots,j_q): i_2,\cdots i_p,j_1,\cdots,j_q\in N\},\\
&&\Delta^S:=\{(i_2,\cdots i_p,j_1,\cdots,j_q): i_2,\cdots i_p,j_1,\cdots,j_q\in S\},\\
&&\Omega^N:=\{(i_1,\cdots i_p,j_2,\cdots,j_q): i_1,\cdots i_p,j_2,\cdots,j_q\in N\},\\
&&\Omega^S:=\{(i_1,\cdots i_p,j_2,\cdots,j_q): i_1,\cdots i_p,j_2,\cdots,j_q\in S\},
\end{eqnarray*}
and then
\[\overline{\Delta^S}=\Delta^N\backslash\Delta^S,~~\overline{\Omega^S}=\Omega^N\backslash\Omega^S.\]
This implies that for a rectangular tensor $\mathcal{A}=(a_{i_1\cdots i_pj_1\cdots j_q}),$ we have that for $i,j\in S$,
\begin{eqnarray*}
R_i(\mathcal{A})=\sum\limits_{i_2,\cdots,i_p,j_1,\cdots,j_q\in N}|a_{ii_2\cdots i_pj_1\cdots j_q}|=r_i^{\Delta^S}(\mathcal{A})+r_i^{\overline{\Delta^S}}(\mathcal{A}),\\
C_j(\mathcal{A})=\sum\limits_{i_1,\cdots,i_p,j_2,\cdots,j_q\in N}|a_{i_1\cdots i_pjj_2\cdots j_q}|=c_j^{\Omega^S}(\mathcal{A})+c_j^{\overline{\Omega^S}}(\mathcal{A}),
\end{eqnarray*}
where
\begin{eqnarray*}
&&r_i^{\Delta^S}(\mathcal{A})=\sum\limits_{(i_2,\cdots,i_p,j_1,\cdots,j_q)\in \Delta^S}|a_{ii_2\cdots i_pj_1\cdots j_q}|,~~r_i^{\overline{\Delta^S}}(\mathcal{A})=\sum\limits_{(i_2,\cdots,i_p,j_1,\cdots,j_q)\in \overline{\Delta^S}}|a_{ii_2\cdots i_pj_1\cdots j_q}|,\\
&&c_j^{\Omega^S}(\mathcal{A})=\sum\limits_{(i_1,\cdots,i_p,j_2,\ldots,j_q)\in \Omega^S}|a_{i_1\cdots i_pjj_2\cdots j_q}|,~~c_j^{\overline{\Omega^S}}(\mathcal{A})=\sum\limits_{(i_1,\cdots,i_p,j_2,\cdots,j_q)\in \overline{\Omega^S}}|a_{i_1\cdots i_pjj_2\cdots j_q}|.
\end{eqnarray*}

\begin{thm}\label{th1}
Let $\mathcal{A}\in \mathbb{R}^{[p,q;n,n]}$, $S$ be a nonempty proper subset of $N$, $\bar{S}$ be the complement of $S$ in $N$. Then
\begin{eqnarray*}
\sigma(\mathcal{A})\subseteq \Psi^S(\mathcal{A})=
\left(\bigcup\limits_{i\in S, j\in \bar{S}}\left(\hat{\Psi}_{i,j}^S(\mathcal{A})\bigcup\tilde{\Psi}_{i,j}^S(\mathcal{A})\right)\right)\bigcup
\left(\bigcup\limits_{i\in \bar{S}, j\in S}\left(\hat{\Psi}_{i,j}^{\bar{S}}(\mathcal{A})\bigcup\tilde{\Psi}_{i,j}^{\bar{S}}(\mathcal{A})\right)\right),
\end{eqnarray*}
where
\begin{eqnarray*}
\hat{\Psi}_{i,j}^S(\mathcal{A})=\left\{z\in \mathbb{C}: (|z|-r_j^{\overline{\Delta^S}}(\mathcal{A}))|z|\leq r_j^{\Delta^S}(\mathcal{A})\max\{R_i(\mathcal{A}),C_i(\mathcal{A})\}\right\},\\
\hat{\Psi}_{i,j}^{\bar{S}}(\mathcal{A})=\left\{z\in \mathbb{C}: (|z|-r_j^{\overline{\Delta^{\bar{S}}}}(\mathcal{A}))|z|\leq r_j^{\Delta^{\bar{S}}}(\mathcal{A})\max\{R_i(\mathcal{A}),C_i(\mathcal{A})\}\right\},\\
\tilde{\Psi}_{i,j}^S(\mathcal{A})=\left\{z\in \mathbb{C}: (|z|-c_j^{\overline{\Delta^S}}(\mathcal{A}))|z|\leq c_j^{\Delta^S}(\mathcal{A})\max\{R_i(\mathcal{A}),C_i(\mathcal{A})\}\right\},\\
\tilde{\Psi}_{i,j}^{\bar{S}}(\mathcal{A})=\left\{z\in \mathbb{C}: (|z|-c_j^{\overline{\Delta^{\bar{S}}}}(\mathcal{A}))|z|\leq c_j^{\Delta^{\bar{S}}}(\mathcal{A})\max\{R_i(\mathcal{A}),C_i(\mathcal{A})\}\right\}.
\end{eqnarray*}
\end{thm}
\noindent {\bf Proof.}
For any $\lambda\in \sigma(\mathcal{A})$, let $x=(x_1,x_2,\cdots,x_m)^\mathrm{T}\in \mathbb{C}^m\backslash\{0\}$ and $y=(y_1,y_2,\cdots,y_n)^\mathrm{T}\in \mathbb{C}^n\backslash\{0\}$ be the left and right associated eigenvectors, that is,
\begin{numcases}{}
\mathcal{A}x^{p-1}y^q=\lambda x^{[l-1]},\label{th1-equ1}\\
\mathcal{A}x^py^{q-1}=\lambda y^{[l-1]}.\label{th1-equ2}
\end{numcases}
Let
\begin{eqnarray*}
&&\hspace*{2cm}|x_t|=\max\{|x_i|:i\in S\},~|x_h|=\max\{|x_i|:i\in \bar{S}\};\\
&&\hspace*{2cm}|y_f|=\max\{|y_i|:i\in S\},~|y_g|=\max\{|y_i|:i\in \bar{S}\};\\
&&w_i=\max\{|x_i|,|y_i|\},i\in N,~w_S=\max\{w_i:i\in S\},~w_{\bar{S}}=\max\{w_i:i\in \bar{S}\}.
\end{eqnarray*}
Then, at least one of $|x_t|$ and $|x_h|$ is nonzero, and at least of $|y_f|$ and $|y_g|$ is nonzero.
We divide the proof into four parts.

Case I: Suppose that $w_S=|x_t|, w_{\bar{S}}=|x_h|$, then $|x_t|\geq |y_t|, |x_h|\geq |y_h|$.

(i) If $|x_h|\geq |x_t|$, then $|x_h|=\max\{w_i:i\in N\}$.
By the $h$th equality in (\ref{th1-equ1}), i.e.,
\begin{eqnarray*}
\lambda x_h^{l-1}&=&\sum\limits_{(i_2,\cdots, i_p,j_1,\cdots, j_q)\in \Delta^S}a_{hi_2\cdots i_pj_1\cdots j_q}x_{i_2}\cdots x_{i_p}y_{j_1}\cdots y_{j_q}\\
&&+\sum\limits_{(i_2,\cdots, i_p,j_1,\cdots, j_q)\in \overline{\Delta^S} }a_{hi_2\cdots i_pj_1\cdots j_q}x_{i_2}\cdots x_{i_p}y_{j_1}\cdots y_{j_q},
\end{eqnarray*}
we have
\begin{eqnarray*}\nonumber
|\lambda| |x_h|^{l-1}&\leq&\sum\limits_{(i_2,\cdots, i_p,j_1,\cdots, j_q)\in \Delta^S}|a_{hi_2\cdots i_pj_1\cdots j_q}||x_{i_2}|\cdots |x_{i_p}||y_{j_1}|\cdots |y_{j_q}|\\\nonumber
&&+\sum\limits_{(i_2,\cdots, i_p,j_1,\cdots, j_q)\in \overline{\Delta^S} }|a_{hi_2\cdots i_pj_1\cdots j_q}||x_{i_2}|\cdots |x_{i_p}||y_{j_1}|\cdots |y_{j_q}|\\\nonumber
&\leq&\sum\limits_{(i_2,\cdots, i_p,j_1,\cdots, j_q)\in \Delta^S}|a_{hi_2\cdots i_pj_1\cdots j_q}||x_t|^{l-1}
+\sum\limits_{(i_2,\cdots, i_p,j_1,\cdots, j_q)\in \overline{\Delta^S}}|a_{hi_2\cdots i_pj_1\cdots j_q}||x_h|^{l-1}\\\nonumber
&=&r_h^{\Delta^S}(\mathcal{A})|x_t|^{l-1}+r_h^{\overline{\Delta^S}}(\mathcal{A})|x_h|^{l-1},
\end{eqnarray*}
i.e.,
\begin{eqnarray}\label{th1-equ3}
(|\lambda|-r_h^{\overline{\Delta^S}}(\mathcal{A}))|x_h|^{l-1}\leq r_h^{\Delta^S}(\mathcal{A})|x_t|^{l-1}.
\end{eqnarray}
If $|x_t|=0$, then $|\lambda|-r_h^{\overline{\Delta^S}}(\mathcal{A})\leq 0$ as $|x_h|>0,$ and it is obvious that
\begin{eqnarray*}
(|\lambda|-r_h^{\overline{\Delta^S}}(\mathcal{A}))|\lambda|\leq 0\leq r_h^{\Delta^S}(\mathcal{A})R_t(\mathcal{A}),
\end{eqnarray*}
which implies that $\lambda\in \hat{\Psi}_{t,h}(\mathcal{A})$.

Otherwise, $|x_t|>0$. From the $t$th equality in (\ref{th1-equ1}), we have
\begin{eqnarray}\nonumber
|\lambda||x_t|^{l-1}&\leq&\sum\limits_{i_2,\cdots,i_p,j_1,\cdots,j_q\in N}|a_{ti_2\cdots i_pj_1\cdots j_q}||x_{i_2}|\cdots |x_{i_p}||y_{j_1}|\cdots |y_{j_q}|\\\nonumber
&\leq&\sum\limits_{i_2,\cdots, i_p,j_1,\cdots,j_q\in N}|a_{ti_2\cdots i_pj_1\cdots j_q}||x_h|^{l-1}\\\nonumber
&=&R_t(\mathcal{A})|x_h|^{l-1},
\end{eqnarray}
i.e.,
\begin{eqnarray}\label{th1-equ4}
|\lambda||x_t|^{l-1}\leq R_t(\mathcal{A})|x_h|^{l-1}.
\end{eqnarray}
Multiplying (\ref{th1-equ3}) by (\ref{th1-equ4}) and noting that $|x_t|^{l-1}|x_h|^{l-1}>0$, we have
\begin{eqnarray*}
(|\lambda|-r_h^{\overline{\Delta^S}}(\mathcal{A}))|\lambda|\leq r_h^{\Delta^S}(\mathcal{A})R_t(\mathcal{A}),
\end{eqnarray*}
which also implies that
$\lambda\in \hat{\Psi}_{t,h}^S(\mathcal{A})\subseteq\bigcup\limits_{i\in S,j\in \bar{S}}\hat{\Psi}_{i,j}^S(\mathcal{A})$.

(ii) If $|x_t|\geq |x_h|$, then $|x_t|=\max\{w_i:i\in N\}$. Similar to the proof of (i), we have
\begin{eqnarray*}
(|\lambda|-r_t^{\overline{\Delta^{\bar{S}}}}(\mathcal{A}))|\lambda|\leq r_t^{\Delta^{\bar{S}}}(\mathcal{A})R_h(\mathcal{A}),
\end{eqnarray*}
which implies that
$\lambda\in \hat{\Psi}_{h,t}^{\bar{S}}(\mathcal{A})\subseteq\bigcup\limits_{i\in \bar{S},j\in S}\hat{\Psi}_{i,j}^{\bar{S}}(\mathcal{A})$.

Case II: Suppose that $w_S=|y_f|, w_{\bar{S}}=|y_g|$, then $|y_f|\geq |x_f|, |y_g|\geq |x_g|$.
If $|y_g|\geq |y_f|$, then $|y_g|=\max\{w_i:i\in N\}$. Similar to the proof of (i) in Case I, we have
\begin{eqnarray*}
(|\lambda|-c_g^{\overline{\Omega^{S}}}(\mathcal{A}))|\lambda|\leq c_g^{\Omega^{S}}(\mathcal{A})C_f(\mathcal{A}),
\end{eqnarray*}
which implies that $\lambda\in \tilde{\Psi}_{f,g}^S(\mathcal{A})\subseteq\bigcup\limits_{i\in S,j\in \bar{S}}\tilde{\Psi}_{i,j}^S(\mathcal{A})$.

If $|y_f|\geq |y_g|$, then $|y_f|=\max\{w_i:i\in N\}$. Similar to the proof of (ii) in Case I, we have
\begin{eqnarray*}
(|\lambda|-c_f^{\overline{\Omega^{\bar{S}}}}(\mathcal{A}))|\lambda|\leq c_f^{\Omega^{\bar{S}}}(\mathcal{A})C_g(\mathcal{A})
\end{eqnarray*}
and
$\lambda\in \tilde{\Psi}_{g,f}^{\bar{S}}(\mathcal{A})\subseteq\bigcup\limits_{i\in \bar{S},j\in S}\tilde{\Psi}_{i,j}^{\bar{S}}(\mathcal{A})$.

Case III: Suppose that $w_S=|x_t|, w_{\bar{S}}=|y_g|$, then $|x_t|\geq |y_t|, |y_g|\geq |x_g|$.
If $|y_g|\geq |x_t|$, then $|y_g|=\max\{w_i:i\in N\}$.
Similarly, we have
\begin{eqnarray*}\label{th1-equ11}
(|\lambda|-c_g^{\overline{\Omega^S}}(\mathcal{A}))\leq c_g^{\Omega^S}(\mathcal{A})R_t(\mathcal{A}),
\end{eqnarray*}
which implies that $\lambda\in \tilde{\Psi}_{t,g}^S(\mathcal{A})\subseteq\bigcup\limits_{i\in S,j\in \bar{S}}\tilde{\Psi}_{i,j}^S(\mathcal{A})$.

If $|x_t|\geq |y_g|$, then $|x_t|=\max\{w_i:i\in N\}$. Similarly, we have
\begin{eqnarray*}\label{th1-equ14}
(|\lambda|-r_t^{\overline{\Delta^{\bar{S}}}}(\mathcal{A}))|\lambda|\leq r_t^{\Delta^{\bar{S}}}(\mathcal{A})C_g(\mathcal{A}),
\end{eqnarray*}
which implies that $\lambda\in \hat{\Psi}_{g,t}^{\bar{S}}(\mathcal{A})\subseteq\bigcup\limits_{i\in \bar{S},j\in S}\hat{\Psi}_{i,j}^{\bar{S}}(\mathcal{A})$.

Case IV: Suppose that $w_S=|y_f|, w_{\bar{S}}=|x_h|$, then $|y_f|\geq |x_f|, |x_h|\geq |y_h|$.
If $|x_h|\geq |y_f|$, then $|x_h|=\max\{w_i:i\in N\}$.
Similarly, we have
\begin{eqnarray*}\label{th1-equ17}
(|\lambda|-r_h^{\overline{\Delta^S}}(\mathcal{A}))|\lambda|\leq r_h^{\Delta^S}(\mathcal{A})C_f(\mathcal{A}),
\end{eqnarray*}
which also implies that
$\lambda\in \hat{\Psi}_{t,h}^S(\mathcal{A})\subseteq\bigcup\limits_{i\in S,j\in \bar{S}}\hat{\Psi}_{i,j}^S(\mathcal{A})$.

If $|y_f|\geq |x_h|$, then $|y_f|=\max\{w_i:i\in N\}$. Similarly, we have
\begin{eqnarray*}\label{th1-equ20}
(|\lambda|-c_f^{\overline{\Omega^{\bar{S}}}}(\mathcal{A}))|\lambda|\leq c_f^{\Omega^{\bar{S}}}(\mathcal{A})R_h(\mathcal{A})
\end{eqnarray*}
and
$\lambda\in \tilde{\Psi}_{h,f}^{\bar{S}}(\mathcal{A})\subseteq\bigcup\limits_{i\in \bar{S},j\in S}\tilde{\Psi}_{i,j}^{\bar{S}}(\mathcal{A})$.

The conclusion follows from Case I, II, III and IV.
\hfill$\blacksquare$

Next, we give the comparison theorem for Theorem \ref{scl-th1} and Theorem \ref{th1}.

\begin{thm}\label{th2}
Let $\mathcal{A}\in \mathbb{R}^{[p,q;n,n]}$, $S$ be a nonempty proper subset of $N$, $\bar{S}$ be the complement of $S$ in N. Then
\[\Psi^S(\mathcal{A})\subseteq \Upsilon^S(\mathcal{A}).\]
\end{thm}
\noindent {\bf Proof.}
Let $z\in \Psi^S(\mathcal{A})$. Then
\begin{eqnarray*}
z\in\bigcup\limits_{i\in S, j\in \bar{S}}\hat{\Psi}_{i,j}^S(\mathcal{A}),~\textmd{or}~z\in\bigcup\limits_{i\in \bar{S}, j\in S}\hat{\Psi}_{i,j}^{\bar{S}}(\mathcal{A}),~\textmd{or}~ z\in\bigcup\limits_{i\in S, j\in \bar{S}}\tilde{\Psi}_{i,j}^S(\mathcal{A}),
~\textmd{or}~z\in\bigcup\limits_{i\in \bar{S}, j\in S}\tilde{\Psi}_{i,j}^{\bar{S}}(\mathcal{A}).
\end{eqnarray*}
Without loss of generality, suppose that $z\in\bigcup\limits_{i\in S, j\in \bar{S}}\hat{\Psi}_{i,j}^S(\mathcal{A})$, i.e.,
there are $u\in S,v\in \bar{S}$, such that
\begin{eqnarray}\label{th3-equ1}
(|z|-r_v^{\overline{\Delta^S}}(\mathcal{A}))|z|\leq r_v^{\Delta^S}(\mathcal{A})\max\{R_u(\mathcal{A}),C_u(\mathcal{A})\}.
\end{eqnarray}
We divide into two cases to prove that $z\in\Upsilon^S(\mathcal{A})$ (for other cases, we can prove them similarly).

Case I. If $r_v^{\Delta^S}(\mathcal{A})\max\{R_u(\mathcal{A}),C_u(\mathcal{A})\}=0$, then
$r_v^{\Delta^S}(\mathcal{A})=0$ or $\max\{R_u(\mathcal{A}),C_u(\mathcal{A})\}=0$.
When $r_v^{\Delta^S}(\mathcal{A})=0$, we have $a_{vu\cdots uu\cdots u}=0$,
$r_v^{\overline{\Delta^S}}(\mathcal{A})=r_v^u(\mathcal{A})$ and
\begin{eqnarray*}
(|z|-r_v^u(\mathcal{A}))|z|&=&(|z|-r_v^{\overline{\Delta^S}}(\mathcal{A}))|z|\leq r_v^{\Delta^S}(\mathcal{A})\max\{R_u(\mathcal{A}),C_u(\mathcal{A})\}=0\\
&=&|a_{vu\cdots uu\cdots u}|\max\{R_u(\mathcal{A}),C_u(\mathcal{A})\},
\end{eqnarray*}
which implies that $z\in\hat{\Upsilon}_{v,u}(\mathcal{A})\subseteq\bigcup\limits_{i\in \bar{S}, j\in S}\hat{\Upsilon}_{i,j}(\mathcal{A})\subseteq\Upsilon^{S}(\mathcal{A}).$

When $\max\{R_u(\mathcal{A}),C_u(\mathcal{A})\}=0$, we have
\begin{eqnarray*}
(|z|-r_v^u(\mathcal{A}))|z|&\leq&(|z|-r_v^{\overline{\Delta^S}}(\mathcal{A}))|z|
\leq r_v^{\Delta^S}(\mathcal{A})\max\{R_u(\mathcal{A}),C_u(\mathcal{A})\}=0\\
&=&|a_{vu\cdots uu\cdots u}|\max\{R_u(\mathcal{A}),C_u(\mathcal{A})\},
\end{eqnarray*}
which also leads to $z\in\hat{\Upsilon}_{v,u}(\mathcal{A})\subseteq\bigcup\limits_{i\in \bar{S}, j\in S}\hat{\Upsilon}_{i,j}(\mathcal{A})\subseteq\Upsilon^{S}(\mathcal{A}).$

Case II. If $r_v^{\Delta^S}(\mathcal{A})\max\{R_u(\mathcal{A}),C_u(\mathcal{A})\}>0$, then by (\ref{th3-equ1}), we have
\begin{eqnarray}\label{th3-equ2}
\frac{|z|-r_v^{\overline{\Delta^S}}(\mathcal{A})}{r_v^{\Delta^S}(\mathcal{A})}\frac{|z|}{\max\{R_u(\mathcal{A}),C_u(\mathcal{A})\}}\leq 1.
\end{eqnarray}
From (\ref{th3-equ2}), we have
\begin{eqnarray}\label{th3-equ3}
\frac{|z|-r_v^{\overline{\Delta^S}}(\mathcal{A})}{r_v^{\Delta^S}(\mathcal{A})}\leq 1
\end{eqnarray}
or
\begin{eqnarray}\label{th3-equ4}
\frac{|z|}{\max\{R_u(\mathcal{A}),C_u(\mathcal{A})\}}\leq 1.
\end{eqnarray}
Let $a=|z|,b=r_v^{\overline{\Delta^S}}(\mathcal{A}),c=r_v^{\Delta^S}(\mathcal{A})-|a_{vu\cdots uu\cdots u}|$ and $d=|a_{vu\cdots uu\cdots u}|$.
When (\ref{th3-equ3}) holds and $d=|a_{vu\cdots uu\cdots u}|>0,$ by Lemma 2.2 in \cite{lcq-lyt} and (\ref{th3-equ2}), we have
\begin{eqnarray*}
\frac{|z|-r_v^u(\mathcal{A})}{|a_{vu\cdots uu\cdots u}|}\frac{|z|}{\max\{R_u(\mathcal{A}),C_u(\mathcal{A})\}}\leq\frac{|z|-r_v^{\overline{\Delta^S}}(\mathcal{A})}{r_v^{\Delta^S}(\mathcal{A})}\frac{|z|}{\max\{R_u(\mathcal{A}),C_u(\mathcal{A})\}}\leq 1,
\end{eqnarray*}
equivalently,
\begin{eqnarray*}
(|z|-r_v^u(\mathcal{A}))|z|\leq |a_{vu\cdots uu\cdots u}|\max\{R_u(\mathcal{A}),C_u(\mathcal{A})\},
\end{eqnarray*}
which implies that $z\in\hat{\Upsilon}_{v,u}(\mathcal{A})\subseteq\bigcup\limits_{i\in \bar{S}, j\in S}\hat{\Upsilon}_{i,j}(\mathcal{A})\subseteq\Upsilon^{S}(\mathcal{A}).$
When (\ref{th3-equ3}) holds and $d=|a_{vu\cdots uu\cdots u}|=0,$ we have
\begin{eqnarray*}
|z|-R_v(\mathcal{A})\leq 0, ~i.e.,~ |z|-r_v^u(\mathcal{A})\leq |a_{vu\cdots uu\cdots u}|=0.
\end{eqnarray*}
Hence,
\begin{eqnarray*}
(|z|-r_v^u(\mathcal{A}))|z|\leq 0=|a_{vu\cdots uu\cdots u}|\max\{R_u(\mathcal{A}),C_u(\mathcal{A})\},
\end{eqnarray*}
which also implies that $z\in\hat{\Upsilon}_{v,u}(\mathcal{A})\subseteq\bigcup\limits_{i\in \bar{S}, j\in S}\hat{\Upsilon}_{i,j}(\mathcal{A})\subseteq\Upsilon^{S}(\mathcal{A}).$

On the other hand, when inequality (\ref{th3-equ4}) holds, i.e., $|z|\leq\max\{R_u(\mathcal{A}),C_u(\mathcal{A})\}$, we only need to prove $z\in\Upsilon^{S}(\mathcal{A})$ under the case that
\begin{eqnarray}\label{th3-equ5}
\frac{|z|-r_v^{\overline{\Delta^S}}(\mathcal{A})}{r_v^{\Delta^S}(\mathcal{A})}> 1,~i.e., \frac{|z|}{R_v(\mathcal{A})}>1.
\end{eqnarray}
When $|z|\leq\max\{R_u(\mathcal{A}),C_u(\mathcal{A})\}=R_u(\mathcal{A})$ and $|a_{uv\cdots vv\cdots v}|>0$,
then from Lemma 2.2, Lemma 2.3 of \cite{lcq-lyt} and (\ref{th3-equ2}), we have
\begin{eqnarray*}
\frac{|z|}{\max\{R_v(\mathcal{A}),C_v(\mathcal{A})\}}\frac{|z|-r_u^v(\mathcal{A})}{|a_{uv\cdots vv\cdots v}|}&\leq&
\frac{|z|}{R_v(\mathcal{A})}\frac{|z|-r_u^v(\mathcal{A})}{|a_{uv\cdots vv\cdots v}|}\\
&\leq&
\frac{|z|-r_v^{\overline{\Delta^S}}(\mathcal{A})}{r_v^{\Delta^S}(\mathcal{A})}\frac{|z|}{R_u(\mathcal{A})}\\
&=&\frac{|z|-r_v^{\overline{\Delta^S}}(\mathcal{A})}{r_v^{\Delta^S}(\mathcal{A})}\frac{|z|}{\max\{R_u(\mathcal{A}),C_u(\mathcal{A})\}}\\
&\leq& 1,
\end{eqnarray*}
equivalently,
\begin{eqnarray*}
(|z|-r_u^v(\mathcal{A}))|z|\leq |a_{uv\cdots vv\cdots v}|\max\{R_v(\mathcal{A}),C_v(\mathcal{A})\},
\end{eqnarray*}
which implies that $z\in\hat{\Upsilon}_{u,v}(\mathcal{A})\subseteq\bigcup\limits_{i\in S, j\in \bar{S}}\hat{\Upsilon}_{i,j}(\mathcal{A})\subseteq\Upsilon^{S}(\mathcal{A}).$
And when $|z|\leq\max\{R_u(\mathcal{A}),C_u(\mathcal{A})\}=R_u(\mathcal{A})$ and $|a_{uv\cdots vv\cdots v}|=0$,
then
$$|z|-r_u^v(\mathcal{A})\leq |a_{uv\cdots vv\cdots v}|=0$$
and
\begin{eqnarray*}
(|z|-r_u^v(\mathcal{A}))|z|\leq 0=|a_{uv\cdots vv\cdots v}|\max\{R_v(\mathcal{A}),C_v(\mathcal{A})\},
\end{eqnarray*}
which also implies that $z\in\hat{\Upsilon}_{u,v}(\mathcal{A})\subseteq\bigcup\limits_{i\in S, j\in \bar{S}}\hat{\Upsilon}_{i,j}(\mathcal{A})\subseteq\Upsilon^{S}(\mathcal{A}).$

When $|z|\leq\max\{R_u(\mathcal{A}),C_u(\mathcal{A})\}=C_u(\mathcal{A})$, similarly, we can obtain that
\begin{eqnarray*}
(|z|-c_u^v(\mathcal{A}))|z|\leq |a_{v\cdots vuv\cdots v}|\max\{R_v(\mathcal{A}),C_v(\mathcal{A})\},
\end{eqnarray*}
which implies that $z\in\tilde{\Upsilon}_{u,v}(\mathcal{A})\subseteq\bigcup\limits_{i\in S, j\in \bar{S}}\tilde{\Upsilon}_{i,j}(\mathcal{A})\subseteq\Upsilon^{S}(\mathcal{A}).$
The proof is completed.
\hfill$\blacksquare$

\begin{remark}\emph{
For a complex tensor $\mathcal{A}\in \mathbb{R}^{[p,q;n,n]}$, the set $\Upsilon^S(\mathcal{A})$ consists of $2|S|(n-|S|)$
sets $\hat{\Upsilon}_{i,j}(\mathcal{A})$ and $2|S|(n-|S|)$ sets $\tilde{\Upsilon}_{i,j}(\mathcal{A})$,
and the set $\Psi^S(\mathcal{A})$ consists of $|S|(n-|S|)$
sets $\hat{\Psi}_{i,j}^S(\mathcal{A})$, $|S|(n-|S|)$
sets $\hat{\Psi}_{i,j}^{\bar{S}}(\mathcal{A})$, $|S|(n-|S|)$
sets $\tilde{\Psi}_{i,j}^S(\mathcal{A})$ and $|S|(n-|S|)$ sets $\tilde{\Psi}_{i,j}^{\bar{S}}(\mathcal{A})$, where $S$ is a nonempty proper subset of $N$. Hence, under the same computations,
$\Psi^S(\mathcal{A})$ captures all singular values of $\mathcal{A}$ more precisely than $\Upsilon^S(\mathcal{A})$.}
\end{remark}


Based on Theorem \ref{th1} and Theorem \ref{th2}, and similar to the proof of Theorem 2 of \cite{scljia} and Theorem 5 of \cite{zjxjias}, the following theorems for the largest singular value of nonnegative rectangular tensors can be obtained easily.

\begin{thm}\label{th3}
Let $\mathcal{A}=(a_{i_{1}\cdots i_{m}})\in \mathbb{R}^{[p,q;n,n]}_{+}$, $S$ be a nonempty proper subset of $N$, $\bar{S}$ be the complement of $S$ in N. Then
\begin{eqnarray*}
\Theta^S(\mathcal{A})\leq\lambda_0\leq \Phi^S(\mathcal{A}),
\end{eqnarray*}
where
\begin{eqnarray*}
\Theta^S(\mathcal{A})=\min\{\hat{\Theta}^S(\mathcal{A}),\hat{\Theta}^{\bar{S}}(\mathcal{A}),\tilde{\Theta}^S(\mathcal{A}),\tilde{\Theta}^{\bar{S}}(\mathcal{A})\}\\ \Phi^S(\mathcal{A})=\max\{\hat{\Phi}^S(\mathcal{A}),\hat{\Phi}^{\bar{S}}(\mathcal{A}),\tilde{\Phi}^S(\mathcal{A}),\tilde{\Phi}^{\bar{S}}(\mathcal{A})\},
\end{eqnarray*}
and
\begin{eqnarray*}
&&\hat{\Phi}^S(\mathcal{A})=\max\limits_{i\in S,j\in \bar{S}}\frac{1}{2}\left\{r_j^{\overline{\Delta^S}}(\mathcal{A})
+[(r_j^{\overline{\Delta^S}}(\mathcal{A}))^2+4r_j^{\Delta^S}(\mathcal{A})\max\{R_i(\mathcal{A}),C_i(\mathcal{A})\}]^{\frac{1}{2}}\right\},\\
&&\hat{\Phi}^{\bar{S}}(\mathcal{A})=\max\limits_{i\in \bar{S},j\in S}\frac{1}{2}\left\{r_j^{\overline{\Delta^{\bar{S}}}}(\mathcal{A})
+[(r_j^{\overline{\Delta^{\bar{S}}}}(\mathcal{A}))^2+4r_j^{\Delta^{\bar{S}}}(\mathcal{A})\max\{R_i(\mathcal{A}),C_i(\mathcal{A})\}]^{\frac{1}{2}}\right\},\\
&&\tilde{\Phi}^S(\mathcal{A})=\max\limits_{i\in S,j\in \bar{S}}\frac{1}{2}\left\{c_j^{\overline{\Delta^S}}(\mathcal{A})
+[(c_j^{\overline{\Delta^S}}(\mathcal{A}))^2+4c_j^{\Delta^S}(\mathcal{A})\max\{R_i(\mathcal{A}),C_i(\mathcal{A})\}]^{\frac{1}{2}}\right\},\\
&&\tilde{\Phi}^{\bar{S}}(\mathcal{A})=\max\limits_{i\in \bar{S},j\in S}\frac{1}{2}\left\{c_j^{\overline{\Delta^{\bar{S}}}}(\mathcal{A})
+[(c_j^{\overline{\Delta^{\bar{S}}}}(\mathcal{A}))^2+4c_j^{\Delta^{\bar{S}}}(\mathcal{A})\max\{R_i(\mathcal{A}),C_i(\mathcal{A})\}]^{\frac{1}{2}}\right\};\\
&&\hat{\Theta}^S(\mathcal{A})=\min\limits_{i\in S,j\in \bar{S}}\frac{1}{2}\left\{r_j^{\overline{\Delta^S}}(\mathcal{A})
+[(r_j^{\overline{\Delta^S}}(\mathcal{A}))^2+4r_j^{\Delta^S}(\mathcal{A})\min\{R_i(\mathcal{A}),C_i(\mathcal{A})\}]^{\frac{1}{2}}\right\},\\
&&\hat{\Theta}^{\bar{S}}(\mathcal{A})=\min\limits_{i\in \bar{S},j\in S}\frac{1}{2}\left\{r_j^{\overline{\Delta^{\bar{S}}}}(\mathcal{A})
+[(r_j^{\overline{\Delta^{\bar{S}}}}(\mathcal{A}))^2+4r_j^{\Delta^{\bar{S}}}(\mathcal{A})\min\{R_i(\mathcal{A}),C_i(\mathcal{A})\}]^{\frac{1}{2}}\right\},\\
&&\tilde{\Theta}^S(\mathcal{A})=\min\limits_{i\in S,j\in \bar{S}}\frac{1}{2}\left\{c_j^{\overline{\Delta^S}}(\mathcal{A})
+[(c_j^{\overline{\Delta^S}}(\mathcal{A}))^2+4c_j^{\Delta^S}(\mathcal{A})\min\{R_i(\mathcal{A}),C_i(\mathcal{A})\}]^{\frac{1}{2}}\right\},\\
&&\tilde{\Theta}^{\bar{S}}(\mathcal{A})=\min\limits_{i\in \bar{S},j\in S}\frac{1}{2}\left\{c_j^{\overline{\Delta^{\bar{S}}}}(\mathcal{A})
+[(c_j^{\overline{\Delta^{\bar{S}}}}(\mathcal{A}))^2+4c_j^{\Delta^{\bar{S}}}(\mathcal{A})\min\{R_i(\mathcal{A}),C_i(\mathcal{A})\}]^{\frac{1}{2}}\right\}.
\end{eqnarray*}
\end{thm}

\begin{thm}\label{th4}
Let $\mathcal{A}=(a_{i_{1}\cdots i_{m}})\in \mathbb{R}^{[p,q;n,n]}_{+}$, $S$ be a nonempty proper subset of $N$, $\bar{S}$ be the complement of $S$ in N. Then
\begin{eqnarray*}
L^S(\mathcal{A})\leq\Theta^S(\mathcal{A})\leq\lambda_0\leq \Phi^S(\mathcal{A})\leq U^S(\mathcal{A}),
\end{eqnarray*}
\end{thm}


In the end, a numerical example is given to verify the theoretical results.

\begin{example}\emph{
Let $\mathcal{A}\in \mathbb{R}^{[2,2;3,3]}_{+}$
with entries defined as follows:
$$a_{1121}=2,a_{1131}=5,a_{1222}=5,a_{2111}=4,a_{2212}=6,$$
$$a_{2232}=1,a_{3111}=8,a_{3222}=9,a_{3313}=8,a_{3323}=10,$$
other $a_{ijkl}=0.$
By computation, we obtain that all different H-singular values of $\mathcal{A}$ are
$-9.9319$, $-9.1950, -9.1673, -7.5415, -2.9723, -2.1787, -1.9725, -1.0403, ~0, ~0.4172, ~3.1373,$\\$3.3683$, $5.7147, 6.9150$ and $18.0755.$
Let $S=\{3\}$ and $\bar{S}=\{1, 2\}$. Next, we consider the singular value inclusion sets and the bounds for the largest singular value of $\mathcal{A}$.}

\emph{(i) $S$-type singular value inclusion sets.}

\emph{By Theorem \ref{scl-th1}, we have
\begin{eqnarray*}
\Upsilon^S(\mathcal{A})=\{z\in{\mathbb{C}}:|z|\leq 33.2547\}.
\end{eqnarray*}
By Theorem \ref{th1}, we have
\begin{eqnarray*}
\Psi^S(\mathcal{A})=\{z\in{\mathbb{C}}:|z|\leq 31.8692\}.
\end{eqnarray*}
The singular value inclusion sets $\Upsilon^S(\mathcal{A})$, $\Psi^S(\mathcal{A})$ and the exact H-singular values are drawn in Figure 1, where
$\Upsilon^S(\mathcal{A})$, $\Psi^S(\mathcal{A})$ and the exact H-singular values are represented by black solid boundary, blue dashed boundary  and red ``$+$", respectively.
It is easy to see that $\Psi^S(\mathcal{A})$ is tighter than $\Upsilon^S(\mathcal{A})$ from Figure 1.}


\begin{figure}[!ht]
\centerline{\includegraphics[height=10cm,width=20cm]{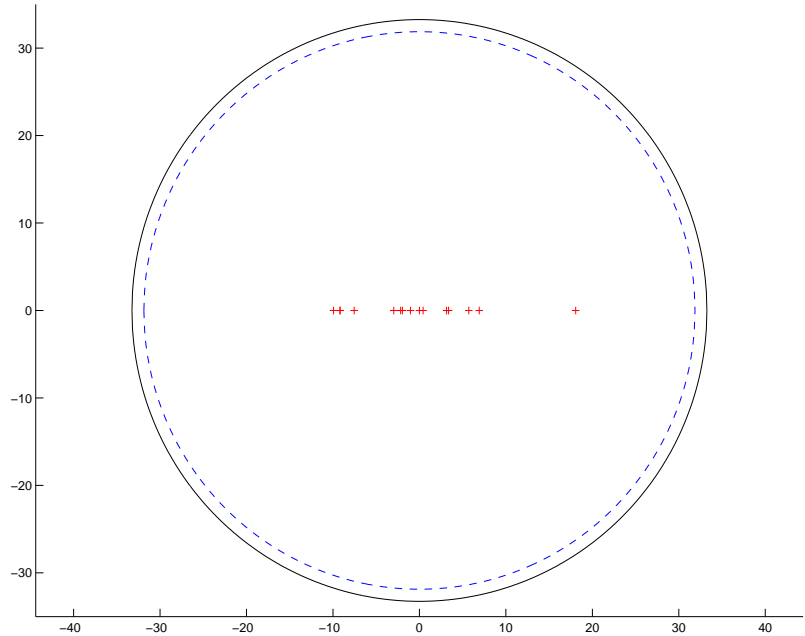}}
{\caption{The comparisons of singular value inclusion sets $\Upsilon^S(\mathcal{A})$ and $\Psi^S(\mathcal{A})$.}}
\label{Fig1}
\end{figure}

\emph{(ii) Bounds for the largest singular value $\lambda_0$.}

\emph{By Theorem 4 of \cite{yangqz}, we have $$6\leq \lambda_0\leq 35.$$
By Theorem \ref{scl-th2}, we have $$6.6533\leq\lambda_0\leq 33.2547.$$
By Theorem \ref{th3}, we have $$8.1240\leq\lambda_0\leq 31.8692.$$
In fact, $\lambda_0=18.0755.$
This example shows that the bounds in Theorem \ref{th3} are better than those in
Theorem \ref{scl-th2} and Theorem 4 of \cite{yangqz}.}
\end{example}
\section{Conclusion}\label{Sec5}

In this paper, by breaking $N=\{1, 2,\cdots,n\}$ into disjoint subsets $S$ and its complement $\bar{S}$, we propose a new $S$-type singular inclusion sets $\Psi^S(\mathcal{A})$ for a real rectangular tensor $\mathcal{A}$ and prove that $\Psi^S(\mathcal{A})$ is tighter than $\Upsilon^S(\mathcal{A})$ in \cite{scljia}.
Based on the set $\Psi^S(\mathcal{A})$, we obtain a new $S$-type upper bound $\Phi^S(\mathcal{A})$ and a new $S$-type lower bound $\Theta^S(\mathcal{A})$ for the largest singular value of nonnegative rectangular tensors, and show that $\Phi^S(\mathcal{A})$ and $\Theta^S(\mathcal{A})$ are sharper than
$U^S(\mathcal{A})$ and $L^S(\mathcal{A})$ in \cite{scljia}, respectively.

Compared with the results in \cite{scljia}, the advantage of the new results is that,
under the same computations, $\Psi^S(\mathcal{A})$ can capture all singular values of $\mathcal{A}$ more precisely than $\Upsilon^S(\mathcal{A})$,
and $\Phi^S(\mathcal{A})$ and $\Theta^S(\mathcal{A})$ can obtain more accurate upper and lower bounds than
$U^S(\mathcal{A})$ and $L^S(\mathcal{A})$, respectively.

\section*{Competing interests}
The author declares that they have no competing interests.
\section*{Authors' contributions}
All authors contributed equally to this work. All authors read and approved the final manuscript.
\section*{Acknowledgments}
This work is supported by National Natural Science Foundations of China (Grant No.11501141),
Foundation of Guizhou Science and Technology Department (Grant No.[2015]2073)
and Natural Science Programs of Education Department of Guizhou Province (Grant No.[2016]066).

\enddocument